\documentclass[12pt,reqno]{amsart}
\usepackage{amsmath,amsfonts,amssymb,amscd,amsthm,amsbsy, color}
\usepackage{appendix}
\usepackage{graphicx}
\usepackage{verbatim}
\usepackage{color}
\usepackage{caption,subcaption}
\usepackage{url}

\usepackage{xcolor}

\usepackage{comment}

\numberwithin{equation}{section}

\theoremstyle{plain}

\theoremstyle{definition}

\newcommand\beq[1]{ \begin{equation}\label{#1} }
\newcommand{\eeq}{ \end{equation} }

\newcommand\beqa[1]{ \begin{eqnarray} \label{#1}}
\newcommand{\eeqa}{ \end{eqnarray} }
\newcommand{\beqano}{ \begin{eqnarray*} }
\newcommand{\eeqano}{ \end{eqnarray*} }

\begin{document}

\title[Periodic orbits near collision]
{Periodic orbits near collision in a restricted four-body problem for the eight choreography}

\author[A. Bengochea]{Abimael Bengochea}
\address{Department of Mathematics, ITAM, R\'io Hondo 1, 01080 Ciudad de M\'exico, M\'exico}
\email{abimael.bengochea@itam.mx}

\author[J. Burgos-Garc\'ia]{Jaime Burgos-Garc\'ia}
\address{Facultad de Ciencias F\'isico Matem\'aticas, Universidad Aut\'onoma de Coahuila, Unidad Campo Redondo, Edificio A, 25020 Saltillo, Coahuila, M\'exico}
\email{jburgos@uadec.edu.mx }

\author[E. P\'erez-Chavela]{Ernesto P\'erez-Chavela}
\address{Department of Mathematics, ITAM, R\'io Hondo 1, 01080 Ciudad de M\'exico, M\'exico}
\email{ernesto.perez@itam.mx}

\baselineskip=18pt              

\begin{abstract}
We study orbits near collision in a non-autonomous restricted planar four-body problem. This restricted problem 
consists of a massless particle moving under the gravitational influence due to three bodies following the figure-eight choreography. We use regularized coordinates in order to deal numerically with motions near collision, and reversing 
symmetries to study both theoretically and numerically certain type of symmetric periodic orbits of the restricted system. The 
symmetric periodic orbits (initial conditions) were computed with the help of some boundary-value problems.
\end{abstract}

\subjclass[2010]{xxxxx, yyyyy, zzzzz}

\maketitle

\tableofcontents

\section{Introduction}

Natural generalizations of the celebrated restricted three-body problem (RTBP), consist of considering two bodies  (primaries) moving in a particular solution of the gravitational two-body problem and adding a massless particle to study its dynamics under the influence of the primaries. It is well known that the general three-body problem admits the classical Euler and Lagrange solutions, where the primaries move either in elliptical or in circular orbits forming collinear and equilateral configurations respectively. One of the first generalizations in the restricted $4$--body problem (three primaries and a massless particle), was the analysis of the dynamic of a massless particle under the gravitational influence of an equilateral configuration, known in the literature as the circular restricted four-body problem, hereafter referred as CRFBP. This problem has revealed to possess interesting dynamics and several works have been performed to study its basic invariant sets such as equilibrium points, periodic orbits, and invariant manifolds. The interested readers can consult 
 \cite{BurgosDelgado}, \cite{BurgosDelgadoII}, \cite{BurgosLessardJames}, \cite{BurgosVI}, \cite{KepleyII}, \cite{Leandro} and references therein. More recently, the so-called eight choreography was found numerically by Moore \cite{Moore}, and its existence was proved by Chenciner and Montgomery \cite{Chenciner-Montgomery}. In contrast with the classical solutions, for this choreography, there is no closed-form formula to describe the orbit.
 
 In \cite{Fujiwara}, the authors investigated the three-body on the lemniscate, they showed that the figure-eight choreography on the lemniscate was not possible considering the Newtonian potential, but it was realized under an inhomogeneous potential. Consequently, there are no special coordinates where the primaries can be considered fixed as in the RTBP or CRFBP. Nevertheless, in \cite{LaraBengochea} the authors considered the three primaries with equal mass moving on the figure-eight solution and a massless particle interacting with this choreography, in this way they define a new restricted four-body problem named \textit{the restricted eight four-body problem}, from now on referred  as REFBP, which defines a non-autonomous dynamical system. In the same work the authors performed a first exploration of symmetric periodic orbits. 
 
Among the different kinds of periodic orbits in a restricted four body problem, two types of them can be understood in terms of a two-body problem. As a first case we consider the test particle located very 
far from the primaries, in such a way that the primaries can be considered as a single body, whose mass is the sum of the masses of the primaries, such that they behave ``almost'' like a two body problem.
The other case is when the test particle is very close to one of the primaries, this is reminiscent of a two-body problem where the dominant interaction with the test particle is due to only one of the primaries. 
This last case is much, much more complicated than the first one, because here, we must regularize the equations of motion to avoid the
problem with the small denominators, which occurs when the massless particle is very close to one of the primaries. This 
is the main subject of this work. 
Compared with the regularization of the circular RTBP and the circular RFBP, the regularization of the current problem is more tricky, because of the lack of suitable coordinates together with the fact that the system is not autonomous. Typically, for the circular and elliptic versions of RTBP and RFBP, the regularization process for collisions between two bodies, relies on the so-called Jacobi integral structures which are, naturally, related to the energy of the system. In our problem such a structure does not exist, nevertheless, in \cite{Broucke} the author shows that a regularization process does not need Jacobi-like integrals or suitable coordinates to be performed. Furthermore, it is enough to consider the \textit{Keplerian} energy between the pair of bodies near collision, this energy has interesting properties and it will be used in our work. Another important aspect of the REFBP is the determination of periodic orbits, since it is more sophisticated than the determination of symmetries. In this way, the isosceles configuration of the primaries on the figure-eight solution will play a main role, as we will show ahead in this paper.
	
We must mention that the study of this problem is twofold; from a mathematical point of view, we believe that the analysis of the periodic orbits in the ERFBP could be used to generate new orbits for the general four-body problem or to study a couple of important open questions: Is it possible to get choreographies with different masses?
Is Saari's conjecture true for more than three bodies? 
From a practical point of view, inspired by recent astronomical observations that reveal that Lagrangian-type configurations can be found beyond the solar system, where Trojan planets are likely to exist, as it was suggested by the observations on the system PDS-70 where a young K7-type star seems to has two giant Jupiter-like protoplanets: PDS 70b and PDS 70c, and since the eight-choreography is KAM stable, it seems likely to be observed in some future time \cite{PDS1},\cite{alma}. Therefore, an available study of particles interacting with this system may result useful.

The paper is organized as follows: in Section \ref{sec:eqmotion} we introduce the equations of motion for the general $4$--body problem, and from here we regularize the equations of motion for the case when the masses $1$ and $4$ are sufficiently close each other, they form a subsystem for the general problem, this allowed us to find new families of symmetric periodic orbits near collision in the REFBP. In Section \ref{sec:Restricted} we formally define the REFBP. In Section \ref{sec:Reversing} we introduce the concepts of symmetry and reversing symmetry of a system, and using a result of Mu\~noz-Almaraz et al. \cite{MunozFreireGalanVader}, which gives a relationship between the fixed points of reversing symmetries and symmetric periodic orbits, we obtain the fixed points associated to a reversing symmetry. Finally, in Section \ref{sec:Numerical}, we do the numerical study for the REFBP for the case when the massless particle is close enough to one of the primaries in the eight configuration. We establish some boundary problems for the non-autonomous system whose sets solution contain initial conditions of symmetric periodic orbits. The numerical solution of boundary problems was obtained with AUTO while high precise initial conditions of periodic orbits were computed with Julia through some root problem. We show several new families of symmetric periodic orbits.

\section{Regularized equations of motion}
\label{sec:eqmotion}
Consider $4$ particles with masses $m_1,...,m_4\geq0$ moving in a two plane where their positions are denoted respectively by the vectors $\textbf{r}_i$, $i=1,...,4$. The equations of motion of the $4$ bodies are
\begin{equation}
\ddot{\textbf{r}}_i=\sum_{j=1, j\neq i}^{4}\frac{m_j(\textbf{r}_j-\textbf{r}_i)}{\vert \textbf{r}_j-\textbf{r}_i\vert^{3}}, \enskip i=1,...,4,
\label{eq:realequations}
\end{equation}
where $\vert\cdot\vert$ is the Euclidian norm of $\mathbb{R}^2$.

The system \eqref{eq:realequations} is not defined if at least two bodies are in the same position, that is collisions. Let
$$
\Delta_{ij} = \left\{({\bf r_1},\cdots,{\bf r}_4) \in \mathbb{R}^8 \enskip | \enskip {\bf r}_i = {\bf r}_j \right\}, \enskip \Delta = \bigcup_{i<j} \Delta_{ij},
$$
then the configuration space of the system is $\Lambda = \mathbb{R}^8 \setminus \Delta$ and $\Omega = \Lambda \times \mathbb{R}^8$ is the corresponding phase space. 

There exist several techniques to regularize collisions in the equations appearing in Celestial Mechanics, the literature dealing with regularizations for restricted problems, where it is possible to construct suitable coordinates such that the primary bodies can be considered fixed, is so vast. Nevertheless, in the case where it is difficult to obtain closed-form formulas for the position of the primaries, as in the problem tackled in this paper, the analysis turns out more interesting, especially if it is desired to use the geometrical advantages of the Hamiltonian formalism. So, it is necessary to extend the odd-dimensional phase space of time-dependent systems to an even-dimensional extended phase space in such a way that the new approach is analogous to the conventional one for autonomous Hamiltonian systems \cite{Celletti}. In the following discussion we will settle the bases of the procedure without a Hamiltonian formalism \cite{Szebehely}. 

We are interested in using suitable coordinates for the study of close approaches between particles 1 and 4. The main idea is to replace ${\bf r}_1$, ${\bf r}_4$ with two planar ${\bf Q}$, ${\bf u}$ vectors defined by
\begin{equation}
{\bf Q}=\frac{m_1 {\bf r}_1 + m_4 {\bf r}_4}{m_1+m_4},\qquad L_u {\bf u} ={\bf r}_4-{\bf r}_1, \qquad
L_u = \left(
\begin{array}{cc}
u_1 & -u_2 \\
u_2 & u_1 \\
\end{array}
\right),
\label{change}
\end{equation}
and incorporate as new coordinate the {\it binding} energy 
\begin{equation}
h = \frac{1}{2}({\bf v}_4-{\bf v}_1)^2 - \frac{m_1+m_4}{|{\bf r}_4-{\bf r}_1|}
\label{hcartesian}
\end{equation}
and fictitious time $\tau$ defined as
\begin{equation}
\frac{d t}{d \tau} = {\bf u}^2.
\label{tau}
\end{equation}
For the new set, instead of differential equations for ${\bf r}_1$ and ${\bf r}_4$, we deal with equations of motion for ${\bf Q}$ and ${\bf u}$ which are obtained using \eqref{eq:realequations}-\eqref{tau}. Notice that the differential equations for ${\bf r}_2$ and ${\bf r}_3$ are keeped, we only replaced the position vectors ${\bf r}_1$, ${\bf r}_4$ by
\begin{equation}
{\bf r}_1 = -\frac{m_4}{m_1+m_4} L_u {\bf u} + {\bf Q}, \qquad {\bf r}_4 = \frac{m_1}{m_1+m_4} L_u {\bf u} + {\bf Q}
\label{r1r4}
\end{equation}
in the right hand side of \eqref{eq:realequations}. In the new set of equations it is also incorporated a new one which describes the temporal evolution of the binding energy. The equations of motion in regularized coordinates are

\begin{equation}
\begin{aligned}
\frac{d{\bf u}}{d  \tau} &= \boldsymbol{\omega}, &\qquad \frac{d {\bf Q}}{d \tau} &= {\bf u}^2 {\bf v}_Q, &\qquad \frac{dh}{d \tau} &= 2 \boldsymbol{\omega}^T L_u^T W, &\qquad \frac{d {\bf r}_i}{d \tau} &= {\bf u}^2 {\bf v}_i, \\[20pt]
 \frac{d \boldsymbol{\omega}}{d  \tau} &= \frac{1}{2}h {\bf u} + \frac{1}{2} {\bf u}^2 L_u^T W, &\qquad \frac{d {\bf v}_Q}{d \tau} &= {\bf u}^2 {\bf a}_Q, &\qquad  \frac{d t}{d \tau} &= {\bf u}^2, &\qquad \frac{d {\bf v}_i}{d \tau} &= {\bf u}^2 {\bf a}_i, 
\end{aligned}
\label{eqmotreg}
\end{equation}
where $i=2,3,$ 
\begin{eqnarray*}
W&=&\sum_{j\neq1,4}m_j\left(\frac{{\bf r}_j-{\bf r}_4}{\vert {\bf r}_j-{\bf r}_4\vert^{3}}-\frac{{\bf r}_j-{\bf r}_1}{\vert {\bf r}_j-{\bf r}_1\vert^{3}}\right),	\\
{\bf a}_Q &=& \sum_{i=1,4,j=2,3}
\frac{m_i m_j}{m_i+m_j} \frac{{\bf r}_j-{\bf r}_i}{\vert {\bf r}_j-{\bf r}_i\vert^{3}},\\
{\bf a}_i &=&\sum_{j \neq i}m_j \frac{{\bf r}_j-{\bf r}_i}{\vert {\bf r}_j-{\bf r}_i\vert^{3}},
\end{eqnarray*}
and ${\bf r}_1$, ${\bf r}_4$ are given by \eqref{r1r4}.

The binding energy in Cartesian coordinates \eqref{hcartesian}, and the one in regularized coordinates
\begin{equation}
h = \frac{2\boldsymbol{\omega}^2 - (m_1+m_4)}{{\bf u}^2},
\label{h-binding}
\end{equation}
are not well defined at binary collision ${\bf u} = {\bf 0}$. Nevertheless, the differential system \eqref{eqmotreg} is regular at collision.

\section{The restricted figure-eight problem}\label{sec:Restricted}

\begin{figure}[h]
\centering
\includegraphics[width=90mm]{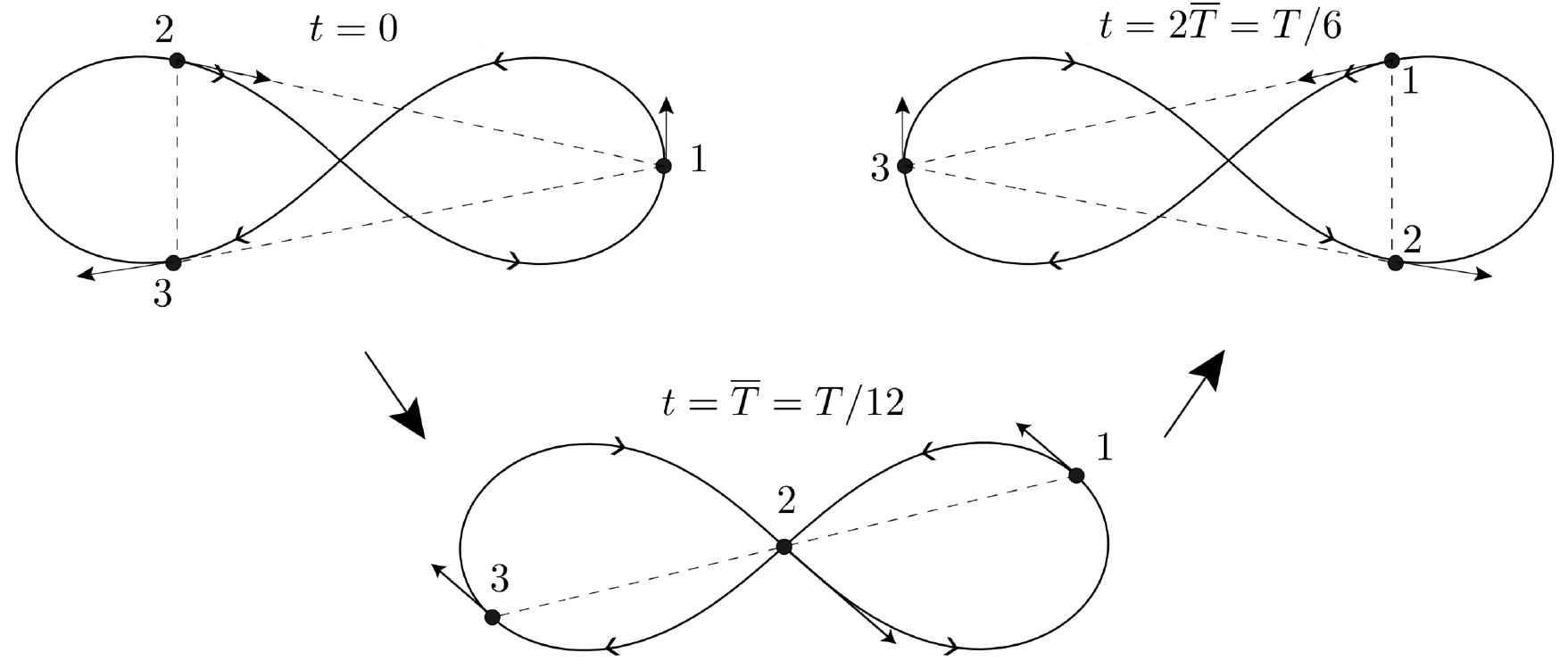}
\caption{Evolution of the figure-eight choreography that starts at isosceles configuration.}
\label{fig1}
\end{figure}

In the following, we describe a non-autonomous restricted four body problem, where the primaries have equal masses and they move in the figure-eight solution, that for short we call it the REFBP  (see \cite{LaraBengochea} for a complete description of this problem). Let us consider four particles on a fixed plane in such a way that the particles $i=1,2,3$, with equal masses (the primaries), follow the figure-eight choreography. The restricted problem consists in the study of the movement of a fourth infinitesimal particle $(m_4 \to 0)$ in presence of the primaries. We use units in such a way that $m_i=1$, $i=1,2,3$, $G=1$, $T=12\overline{T} = 6.32591398$, where $T$ is the period of the choreography, and $2 \overline{T}$ the time that passes between two consecutive isosceles configurations (see Fig.~\ref{fig1}).

In this problem, at time $t=0$ the primary bodies appear in an isosceles configuration, as it is shown in Fig.~\ref{fig1}. According to the choice of units, the positions and velocities of the primaries at the initial time are
\begin{equation}
\begin{array}{l}
{\bf r}_1^T(0) = (1.081017082650648, 0),\\[6pt]
{\bf r}_2^T(0) = (-0.540508541325324, 0.345263314425768),\\[6pt]
{\bf r}_3^T(0) = (-0.540508541325324, -0.345263314425768),\\[6pt]
{\bf v}_1^T(0) = (0, 0.467209527201224),\\[6pt]
{\bf v}_2^T(0) = (1.097122382351121, -0.233604763600612),\\[6pt]
{\bf v}_3^T(0) = (-1.097122382351121, -0.233604763600612).
\end{array}
\label{cieight}
\end{equation}
Later we will require the initial vectors of the first body \eqref{cieight}, so we define its Cartesian entries as ${\bf r}_1^T(0) = (x_{10},0)$, ${\bf v}_1^T(0) = (0,v_{y10})$.

The position and velocity vectors of the primaries ${\bf r}_i (t)$, ${\bf v}_i (t)$, $i = 1, 2, 3$, are known functions through numerical integration of the equations of motion. These vectors are $T$--periodic and independent of the state vectors of the test particle.

The equations of motion of the REFBP (test particle) are
\begin{equation}
\dot{{\bf r}} = {\bf v}, \qquad \dot{\bf v} = \sum_{j=1}^3 \frac{{\bf r}_j(t)-{\bf r}}{||{\bf r}_j(t)-{\bf r}||},
\label{eq-restricted}
\end{equation}
which form a non-autonomous system. In the following we use interchangeably 
$$
{\bf z} =
\left(
\begin{array}{c}
{\bf r} \\
{\bf v}
\end{array}
\right), 
$$
and the cartesian representation
$$
{\bf r} = \left(
\begin{array}{c}
x \\
y
\end{array}
\right), \quad
{\bf v} = \left(
\begin{array}{c}
v_x \\
v_y
\end{array}
\right), \quad
{\bf r}_i = \left(
\begin{array}{c}
x_i \\
y_i
\end{array}
\right), \quad
{\bf v}_i = \left(
\begin{array}{c}
v_{xi} \\
v_{yi}
\end{array}
\right), \quad i=1,2,3.
$$

There is another representation for the REFBP, useful for the numerical implementation: the restricted problem and a four-body problem where the initial conditions of the primaries are given by \eqref{cieight}, and the fourth body has zero mass, of course they are equivalent.

\section{Reversing Symmetries and periodic orbits}\label{sec:Reversing}

For the study of symmetric periodic orbits of the REFBP we use the tool of reversing symmetries and its fixed points \cite{Lamb}. In the following, we introduce relevant concepts and results. 

Consider the system
\begin{equation}
\frac{d{\bf z}}{dt} = F {\bf (z)}, \enskip F: \Omega \longmapsto \Omega,
\label{eq-ust}
\end{equation}
where ${\bf z}$ is the state vector, $F$ is the vector field and $\Omega$ is the phase space. We say that an involution $S: \Omega \longmapsto \Omega$ ($S^2 = Id$), is a symmetry of the system \eqref{eq-ust} if
$$
\frac{dS({\bf z})}{dt} = F \circ S {\bf (z)}
$$
is fulfilled, and a involution $R: \Omega \longmapsto \Omega$ is a reversing symmetry if
$$
\frac{dR({\bf z})}{dt} = -F \circ R{\bf (z)}
$$
holds.

The fixed points of $R$, denoted by
$$
\textnormal{Fix}(R) = \{{\bf z} \in \Omega \enskip \vert \enskip R{\bf z} = {\bf z} \},
$$
are closely related with symmetric periodic solutions as we explain in the following.

We say that a solution ${\bf z}(t)$ of \eqref{eq-ust}, defined in its maximal domain, is $R$--symmetric if its corresponding orbit is invariant under $R$; a well known result is that an orbit is periodic and $R$--symmetric if and only if intersects two points of Fix$(R)$. There are similar concepts and results for two different reversing symmetries $R$ and $\widetilde{R}$. We say that a solution ${\bf z}(t)$ is $(R,\widetilde R)$--symmetric if there exist times $t_0<t_1$ in such a way ${\bf z}(t_0) \in \textnormal{Fix}(R)$ and ${\bf z}(t_1) \in \textnormal{Fix}(\widetilde R)$, where $[t_0, t_1]$ is the basic domain of the doubly-symmetric solution. The following result establishes  the relationship between the fixed points of reversing symmetries and symmetric periodic orbits, you can find a proof of it in \cite{MunozFreireGalanVader}.

The result states that if 
 ${\bf z}(t)$ is a $(R, \widetilde{R})$--symmetric solution with basic domain $[0, T_0]$,
then
\begin{enumerate}
\item ${\bf z}(t)$ is defined for all $t \in \mathbb{R}$,
\item for all $t \in \mathbb{R}$, $m \in Z$, the relations
$$
{\bf z}(-t) = R{\bf z}(t), \enskip {\bf z}(t) = \widetilde{R}{\bf z}(2T_0 - t), \enskip
{\bf z}(2mT_0 + t) = (\widetilde{R} R)^m{\bf z}(t)
$$
are fulfilled.
\end{enumerate}

Notice that if there exists some $M \in \mathbb{N}$ such that $(\widetilde{R} R)^M = id$, the solution is $(R, \widetilde{R})$--symmetric with period $T=2MT_0$.

\subsection{Fixed points of reversing symmetries}

\begin{figure}[h]
\centering
\includegraphics[width=90mm]{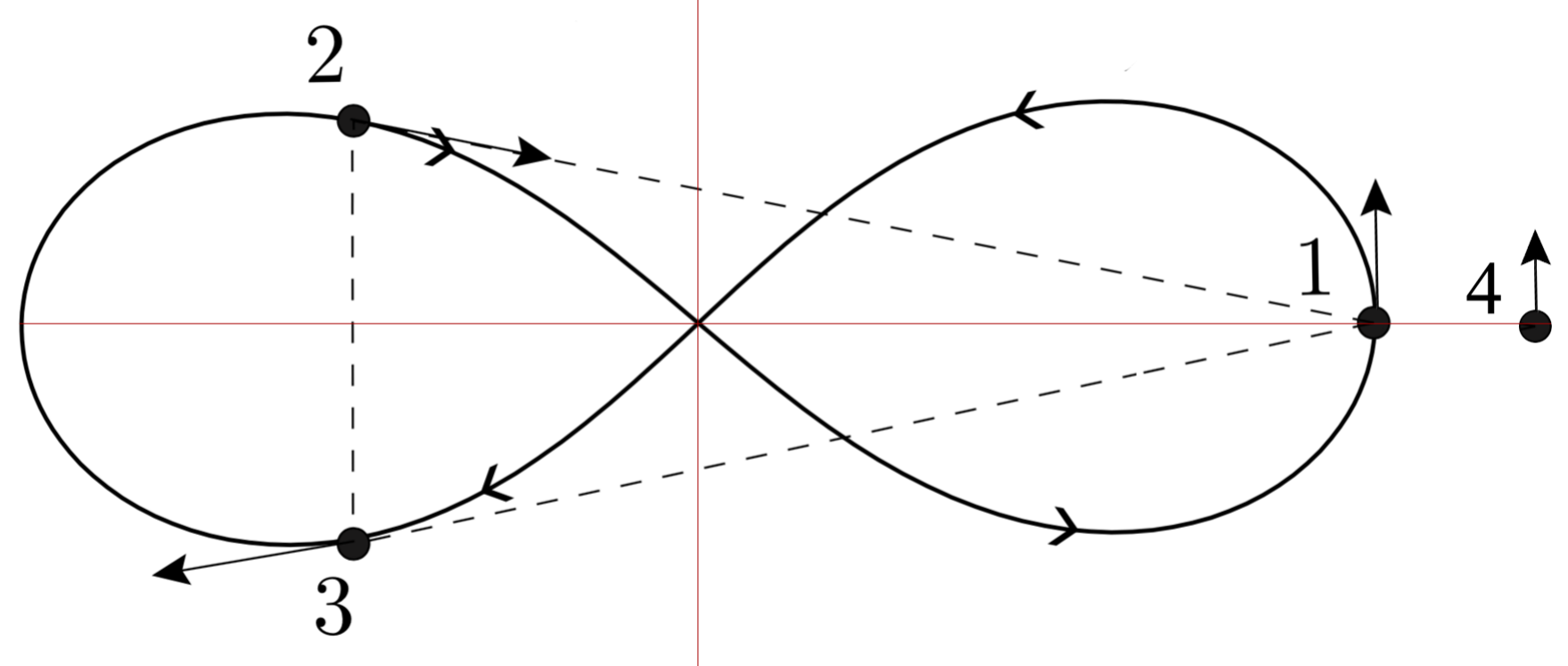}
\caption{Fixed point of reversing symmetry $\widetilde R$ in the REFBP.}
\label{fig2}
\end{figure}

It was proved in \cite{LaraBengochea} that the REFBP possess a reversing symmetry (let us call it $\widetilde R$) in which the primaries are located at an isosceles configuration\footnote{This configuration is a fixed point of a reversing symmetry of the three-body problem.} (see Fig.~\ref{fig1}), and the test particle is on the horizontal axis with velocity parallel to the $y$-axis. According to the description of the REFBP, an orbit could pass through fixed points of $\widetilde R$ only at times $6n \overline{T}$, $n \in \mathbb{Z}$ ($2 \overline{T}$ is the time that passes between two consecutive isosceles configurations of the figure-eight choreography - see Fig.~\ref{fig1}). The equations that describe the fixed points of the reversing symmetry $\widetilde R$ for the test particle are
\begin{equation}
y=0, \enskip v_x=0,
\label{eq-fixedpointrest}
\end{equation}
for $x$ and $v_y$ arbitrary. The corresponding configuration is represented in Fig.~\ref{fig2}

\section{Numerical exploration of symmetric periodic orbits near collision}\label{sec:Numerical}

We are interested in the study of symmetric periodic orbits where the test particle is always near the first primary body, we mean in a moon type periodic orbit. The fixed point of the reversing symmetry $\widetilde R$ are consistent with this type of motion; a solution ${\bf u}(t)$ can only pass through fixed points of $\widetilde R$ at times $t = 6n\overline{T}$, $t \in \mathbb{Z}$. For our analysis, we consider only orbits in which the temporal separation between the corresponding fixed points of $\widetilde R$ is minimal; without loss of generality we choose those associated to times $t=0$ and $t=6 \overline T$. According to \eqref{eq-fixedpointrest}, at these times the test particle lies along the horizontal axis with velocity in the $y$-axis direction, therefore the conditions
\begin{equation}
y(0) = 0, \enskip v_x(0) = 0, \enskip y(6\overline{T}) = 0, \enskip v_x(6\overline{T}) = 0,
\label{t0F}
\end{equation}
must be fulfilled. In terms of the flow $\phi$ of the restricted problem \eqref{eq-restricted}, the system \eqref{t0F} can be stated as follows: find unknowns $x_{0},v_{y0} \in \mathbb{R}$ in such a way that
\begin{equation}
\phi_{6\overline{T}}(x_{0},0,0,v_{y0}) = (x_{f},0,0,v_{yf})
\label{numflow}
\end{equation}
holds, where $x_{f},v_{yf}$ are arbitrary real numbers. Notice that $x_{0},v_{y0}$ represents the initial condition in the $x$--coordinate and $y$--velocity of the test particle, that is $x(0) = x_{0}$, $v_y(0) = v_{y0}$, and $x(6\overline{T}) = x_{f}$, $v_y(6\overline{T}) = v_{yf}$. According to \cite{MunozFreireGalanVader}, a pair $x_{0}$,$v_{y0}$ which satisfies the result of the previous section leads to a symmetric periodic orbit with period $T=12\overline{T}$ \eqref{numflow}.

\begin{figure}[h]
\centering
\includegraphics[width=75mm]{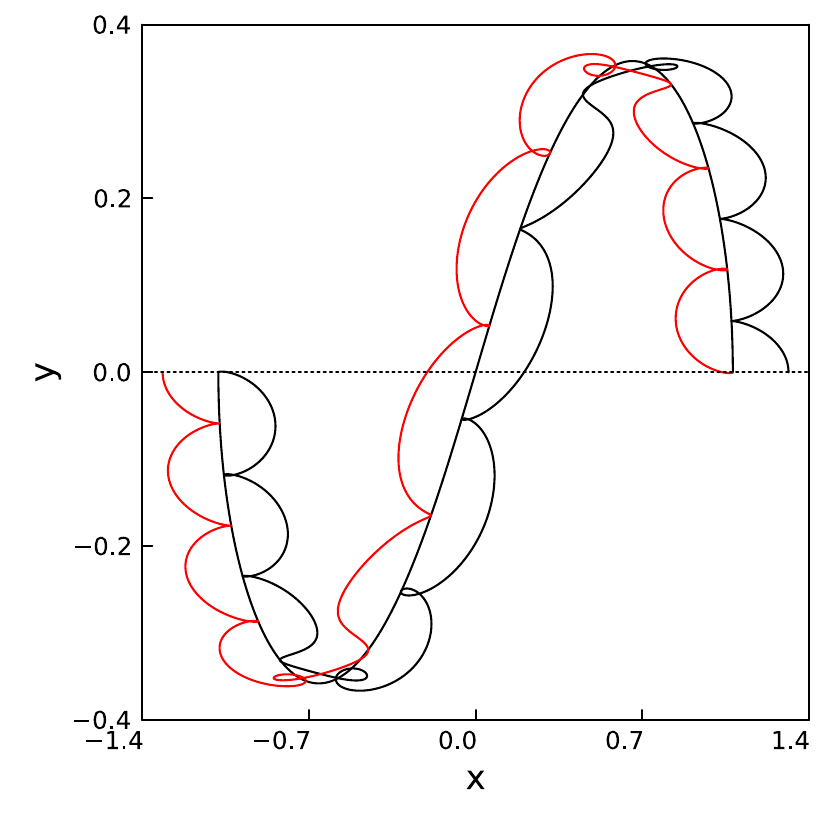}
\caption{Orbits of the REFBP in a half period. A half eight (primary's orbit), and orbits of the test particle with initial conditions $(x_{0},0,0,v_{y0})$, $(-x_{f},0,0,v_{yf})$ are shown (black and red trajectories, respectively).}
\label{figorb1}
\end{figure}

By symmetry, if $(x_{0},0,0,v_{y0})$ satisfies \eqref{numflow} also $(-x_{f},0,0,v_{yf})$ does it. Although these points are different initial conditions of the restricted problem, they represent the same physical orbit, since the initial conditions are related by a reflection along the $y$-axis, and a permutation of the primaries (see Fig. \ref{figorb1}).

\subsection{Boundary Value Problems (Cartesian coordinates)}
\label{BVPcartesian}

With the aim of determining orbits that satisfy \eqref{numflow}, in the following we set three boundary value problems (BVPs) whose set of solutions are denoted by $C_y$, $C_{vx}$ and $C_R$. The motivation to work with these sets is that the pairs $x_{0}$, $v_{y0}$ which belong to the intersection $C_y \cap C_{vx}$ fulfill \eqref{numflow}. In a similar way, $C_R$ is conformed by triads $x_{0}$, $v_{y0}$, $T_0$ that satisfy \eqref{numflow} for specific values of $T_0$ (characteristic time). In the following we describe the BVPs and their solutions in the cartesian coordinate system, in the next section we introduce the equations that define the BVPs in regularized coordinates.

Let us consider the REFBP \eqref{eq-restricted}. We define as first BVP the one in which the solution ${\bf u}(t)$ meets three conditions of \eqref{t0F}, namely
\begin{equation}
\begin{array}{cc}
y(0) = 0, & y(6{\overline T}) = 0, \\ [6pt]
v_{x}(0) = 0.
\end{array}
\label{bvp-1}
\end{equation}
The solution's set of \eqref{bvp-1} is
\begin{equation}
C_y = \left \{ (x_{0},v_{y0}) \in \mathbb{R}^2\, | \, \phi_{6 \overline T}(x_{0},0,0,v_{y0}) = ((x_{f},0,v_{xf},v_{yf})) \right \}
\label{cy}
\end{equation}
with $x_{f}$, $v_{xf}$, $v_{yf}$ arbitrary real numbers.

For the second BVP, the solution ${\bf u}(t)$ of \eqref{eq-restricted} also must satisfy three conditions of \eqref{t0F} given by
\begin{equation}
\begin{array}{c}
\begin{array}{cc}
y(0) = 0, & v_{x}(6{\overline T}) = 0, \\ [6pt]
v_{x}(0) = 0.
\end{array}
\end{array}
\label{bvp-2}
\end{equation}
The solution's set of \eqref{bvp-2} is
\begin{equation}
C_{vx} = \left \{ (x_{0},v_{y0}) \in \mathbb{R}^2\, | \, \phi_{6 \overline T}(x_{0},0,0,v_{y0}) = (x_{f},y_{f},0,v_{yf}) \right \}
\label{cvx}
\end{equation}
where $x_{f}$, $y_{f}$, $v_{yf}$ are arbitrary real numbers.

In general \eqref{cy} and \eqref{cvx} are conformed by curves in the plane $x_{0}v_{y0}$. Moreover, if $(x_{0},v_{y0}) \in C_y \cap C_{v_x}$ then $(x_{0},0,0,v_{y0})$ satisfies \eqref{numflow} and it is an initial condition that leads to a symmetric periodic orbit with period $T=12\overline{T}$.

There is a third BPV, from which we also could obtain points that give rise to initial conditions that satisfy \eqref{numflow}. This BPV is defined as the one in which the solution ${\bf u}(t)$ of \eqref{eq-restricted} satisfies
\begin{equation}
\begin{array}{cc}
y(0) = 0, & y(T_0) = 0,  \\ [6pt]
v_{x}(0) = 0, & v_{x}(T_0) = 0,
\end{array}
\label{bvp-3}
\end{equation}
where $T_0$ is an arbitrary positive real number. Although the only restriction about $T_0$ is that it must be a positive number, in our problem we are dealing only with numerical values near $6\overline T$. The solution's set of \eqref{bvp-3} is given by
\begin{equation}
C_{R} = \left \{  (x_{0},v_{y0},T_0) \in \mathbb{R}^2 \times \mathbb{R}_{>0} \, | \, \phi_{T_0}(x_{0},0,0,v_{y0}) = (x_{f},0,0,v_{yf}) \right \},
\label{cR}
\end{equation}
with arbitrary real numbers $x_{f},v_{yf}$ and $T_0$ positive real number. Notice that if $(x_{0},v_{y0},T_0)$ is a point within $C_{R}$ such that $T_0 = 6 \overline T$ then $(x_{0},0,0,v_{y0})$ is an initial condition of a symmetric periodic orbit since it satisfies \eqref{numflow}.

\subsection{Boundary Value Problems (regularized coordinates)}
\label{BVPregularized}

The equations of the first BVP \eqref{bvp-1}, in function of the regularized coordinates become
\begin{equation}
\begin{array}{cc}
\displaystyle u_1(0)u_2(0) = 0, & u_1(\tau_0)u_2(\tau_0) = 0, \\ [6pt]
\displaystyle u_1(0)\omega_1(0)-u_2(0)\omega_2(0) = 0, & \int_0^{\tau_0} \boldsymbol{u}^2 d\tau = 6 \overline{T}. 
\end{array}
\label{bvp-1-new}
\end{equation}
Notice the correspondence between the first column and the first equation of the second column of \eqref{bvp-1} and \eqref{bvp-1-new}. The second equation in the second column of \eqref{bvp-1-new}, corresponds to the time condition $t=6\overline{T}$. The analogue set to $C_y$, denoted by $\widehat C_y$, is the one made up of points $(u_{10},u_{20},\omega_{10},\omega_{20},\tau_0)$ such that \eqref{bvp-1-new} is fulfilled.

For the second BVP \eqref{bvp-2} the corresponding equations are given by the expressions
\begin{equation}
\begin{array}{cc}
\displaystyle u_1(0)u_2(0) = 0, & u_1(\tau_0)\omega_1(\tau_0)-u_2(\tau_0)\omega_2(\tau_0) = 0, \\ [6pt]
\displaystyle u_1(0)\omega_1(0)-u_2(0)\omega_2(0) = 0, & \int_0^{\tau_0} {\bf u}^2 d\tau = 6\overline{T}.
\end{array}
\label{bvp-2-new}
\end{equation}
Here, the correspondence is between the first column and the first equation of the second column of equations \eqref{bvp-2} and \eqref{bvp-2-new}. As the first BVP, the second equation in second column of \eqref{bvp-2-new} is equivalent to $t=6\overline{T}$. In a similar way $\widehat C_{vx}$, the analogue set of $C_{vx}$, is made up of points $(u_{10},u_{20},\omega_{10},\omega_{20},\tau_0)$ such that \eqref{bvp-2-new} is fulfilled.

The third BVP \eqref{bvp-3} is defined by the equations 
\begin{equation}
\begin{array}{ll}
u_1(0)u_2(0) = 0, & 2u_1(\tau_0)u_2(\tau_0) + y_1(\tau_0)= 0, \\ [6pt]
u_1(0)\omega_1(0)-u_2(0)\omega_2(0) = 0, & 2(u_1(\tau_0)\omega_1(\tau_0)-u_2(\tau_0)\omega_2(\tau_0)) + u^2(\tau_0)v_{x1}(\tau_0) = 0.
\end{array}
\label{bvp-3-new}
\end{equation}
Here, the correspondence holds for each equation of \eqref{bvp-3} and \eqref{bvp-3-new}. The analogue set to $C_{R}$, which we define as $\widehat C_{R}$, is the one made up of points $(u_{10},u_{20},\omega_{10},\omega_{20},\tau_0)$ such that \eqref{bvp-3-new} is fulfilled. Notice that for any vector $(u_{10},u_{20},\omega_{10},\omega_{20},\tau_0)$ within $\widehat C_y$, two of its entries are zero: if $x_0 - x_{10}>0$ then $u_{20} = \omega_{10} = 0$ whereas for $x_0 - x_{10}<0$ we have $u_{10} = \omega_{20} = 0$. Something similar happens for the sets $\widehat C_{vx}$, $\widehat C_{R}$. Thus, these sets are represented in the spaces $u_{10}\omega_{20}\tau_0$ and $u_{20}\omega_{10}\tau_0$.

Just like in the Cartesian case, the initial conditions belonging to $(u_{10},u_{20},\omega_{10},\omega_{20},\tau_0) \in \widehat C_{y} \cap \widehat C_{vx}$ define symmetric periodic orbit with period $12 \overline{T}$. In a similar way, if a point $(u_{10},u_{20},\omega_{10},\omega_{20},\tau_0) \in \widehat C_{R}$ satisfies $\int_0^{\tau_0} {\bf u}^2 d\tau = 6 \overline{T}$ then the initial condition $(u_{10},u_{20},\omega_{10},\omega_{20})$ defines a symmetric periodic orbit with period $12 \overline{T}$.

\subsection{Numerical study}
\label{numstudy}

In order to compute solutions of the BVPs previously defined, that is the $C$ sets either in Cartesian or regularized coordinates, we need a ``seed'' of the corresponding problem. In the following we explain how to proceed in Cartesian coordinates because of the physical interpretation, later we expose the corresponding description in regularized coordinates.

\subsubsection{Seeds for BVPs in Cartesian coordinates}
\label{BVPcartesian}
We consider $C_y$ (see \eqref{cy}). This set is conformed by points $x_{0}$, $v_{y0}$ (initial conditions) such that $y(t)$ is zero at $t = 6 \overline{T}$. Thinking in $y$ also as function of  the initial conditions we have that
\begin{equation}
y(6 \overline{T},x_{0},v_{y0}) = 0
\label{root}
\end{equation}
is equivalent to the condition that defines $C_y$. In \eqref{root}, the parameters $x_{0},v_{y0}$ are unknown; in a generic way $C_y$ is conformed by curves in the plane $x_{0}v_{y0}$.

In order to compute $(x_0,v_{y0}) \in C_y$ we need an approximate value of it, namely $(x_{*},v_{y*})$. For motions near collision we can obtain an approximate point $(x_{*},v_{y*})$ as follows. The main influence over the test particle is due to the first body hence, in the corresponding two body approximation with $m_2 \to 0$, we have elliptic orbits for the bounded case. According to the configuration of initial conditions within $C_y$, the corresponding ellipse must be oriented with its major and minor axes parallel to $x$ and $y$, therefore at the initial time the body 1 and the test particle conform approximately a periapsis/apoapsis configuration and satisfy
\begin{equation}
\begin{array}{l}
\displaystyle {\bf r}_{rel} = a(1 \pm e) \enskip \hat{\imath}, \enskip {\bf v}_{rel} = \sqrt{\frac{1 \mp e}{a(1 \pm e)}} \enskip \hat{\jmath},
\end{array}
\label{rv-rel}
\end{equation}
where $a$ and $e$ are the usual semimajor axis and eccentricity, respectively. In the equation \eqref{rv-rel} we have relative coordinates so we have to add them the corresponding state vectors of the first body to obtain $(x_{*},v_{y*})$.

 Suppose we have obtained $(x_{*},v_{y*})$, if we fix one of the variables, for instance $x_{0} = x_{*}$, we only need to find $v_{y0}$; usually with Newton's method, and a good enough approximation $(x_{*},v_{y*})$, we will be able to compute $(x_{0},v_{y0})$ such that \eqref{root} holds. A similar root problem can be established for obtaining a seed for $C_{vx}$ which also is conformed by curves in the plane $x_{0}v_{y0}$. Once we have a seed for $C_y$ and $C_{vx}$, we can compute some part of the sets. We remember that if $(x_{0},v_{y0}) \in C_{y} \cap C_{vx}$ then such point leads to a symmetric periodic orbit; in particular $(x_{0},v_{y0},6\overline{T})$ is a seed for $C_R$.

\subsubsection{Seeds for BVPs in regularized coordinates}
\label{BVPregularized}

For the study of orbits near collision, we consider the equations of motion in regularized coordinates and the sets $\widehat{C}_{y}$, $\widehat{C}_{vx}$ and $\widehat{C}_{R}$. Each set is represented with vectors $(u_{10},u_{20},\omega_{10},\omega_{20},\tau_0)$ and, as happens in the Cartesian case, is conformed by curves. In all cases the coordinates $u_{10}$, $u_{20}$, $\omega_{10}$, $\omega_{20}$ can be obtained from the associated Cartesian coordinates as we explain in the following.
\begin{enumerate}
\item Case $x_0 - x_{10}>0$. Here immediately follows that  $u_{20} = 0$, $\omega_{10}= 0$. For the initial condition in $u_1$ we have two options $u_{10} = \pm \sqrt{x_{0}-x_{10}}$, and for any election the remaining velocity becomes $\omega_{20} = \frac{1}{2}u_{10}(vy_{0}-vy_{10})$. In this case the initial conditions $(u_{10},0,0,\omega_{20})$ and $(-u_{10},0,0,-\omega_{20})$ are associated to the same orbit in a physically sense.
\item Case $x_0 - x_{10}<0$. For this case, we have $u_{10} = 0$, $\omega_{20} = 0$, and for $u_2$ at the initial time we also have two options, namely $u_{20} = \pm \sqrt{x_{10}-x_{0}}$. As happens in the previous case, for any election we have $\omega_{10} = \frac{1}{2}u_{20}(vy_{0}-vy_{10})$ which lead to initial conditions $(0,u_{20},\omega_{10},0)$ and $(0,-u_{20},-\omega_{10},0)$. Here we also have that these initial conditions define the same orbit in a physical sense.
\end{enumerate}

For the numerical study we always choose the initial condition in such a way that $u_{10}>0$ or $u_{20}>0$. The two possible representations for the regularized sets, namely $u_{10}\omega_{20}$ and $u_{20}\omega_{10}$, are due to the fact that the regularization transformation is not injective.

Once we have the coordinates for the seed $u_{10}$, $u_{20}$, $\omega_{10}$, $\omega_{20}$, we require the characteristic time $\tau_0$. For $\widehat{C}_{y}$ and $\widehat{C}_{vx}$ the value of $\tau_0$ is defined implicitly by the equation $\int_0^{\tau_0} {\bf u}^2 d\tau = 6\overline{T}$, whereas for $\widehat{C}_{R}$ the corresponding condition is $\int_0^{\tau_0} {\bf u}^2 d\tau = T_0$.

\subsubsection{Computing $\widehat{C}_{R}$}

For the purpose of obtaining initial conditions of symmetric periodic orbits is more simple deal with $\widehat C_R$, instead of two sets $\widehat C_y$ and $\widehat C_{vx}$. The main idea that we follow is to determine numerically some curves that conform $\widehat C_R$, and later identify those points belonging $\widehat C_R$ such that $\int_0^{\tau_0} {\bf u}^2 d\tau = 6 \overline{T}$ holds. We remark that there are two type of points within $\widehat C_R$, namely $u_{10}$, $\omega_{20}$, $\tau_0$ or $u_{20}$, $\omega_{10}$, $\tau_0$.

\begin{figure}[h!]
\centering
\includegraphics[width=75mm]{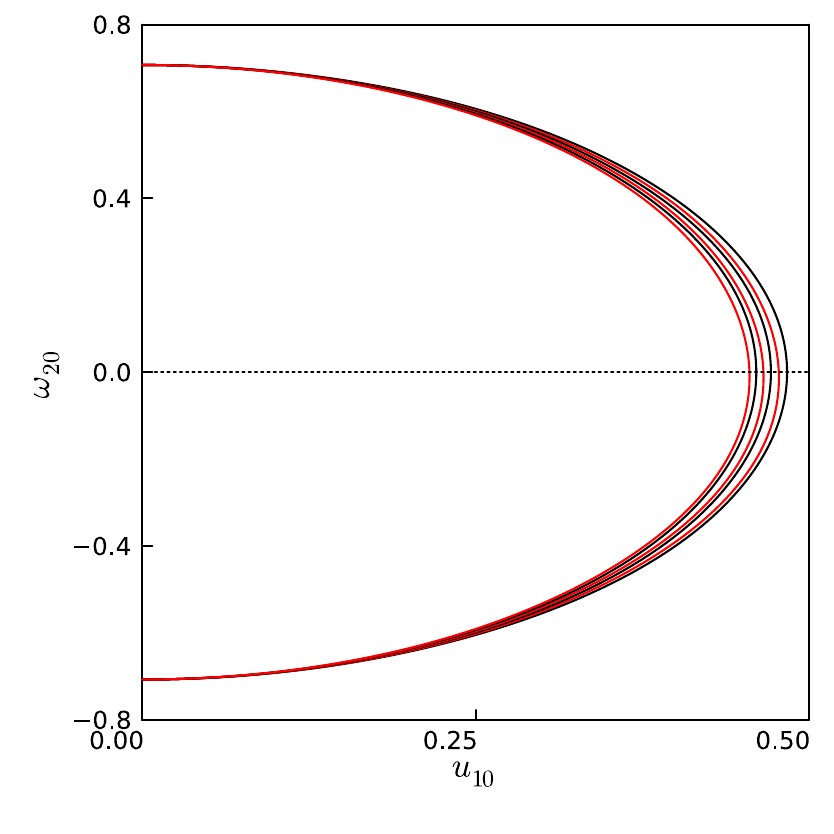}
\includegraphics[width=75mm]{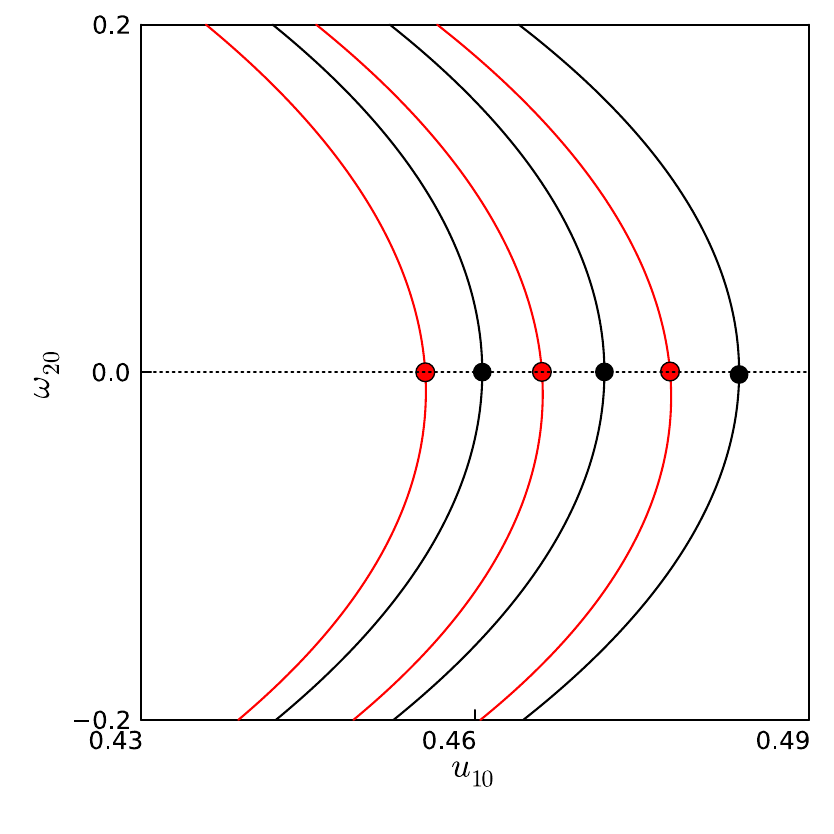}
\includegraphics[width=75mm]{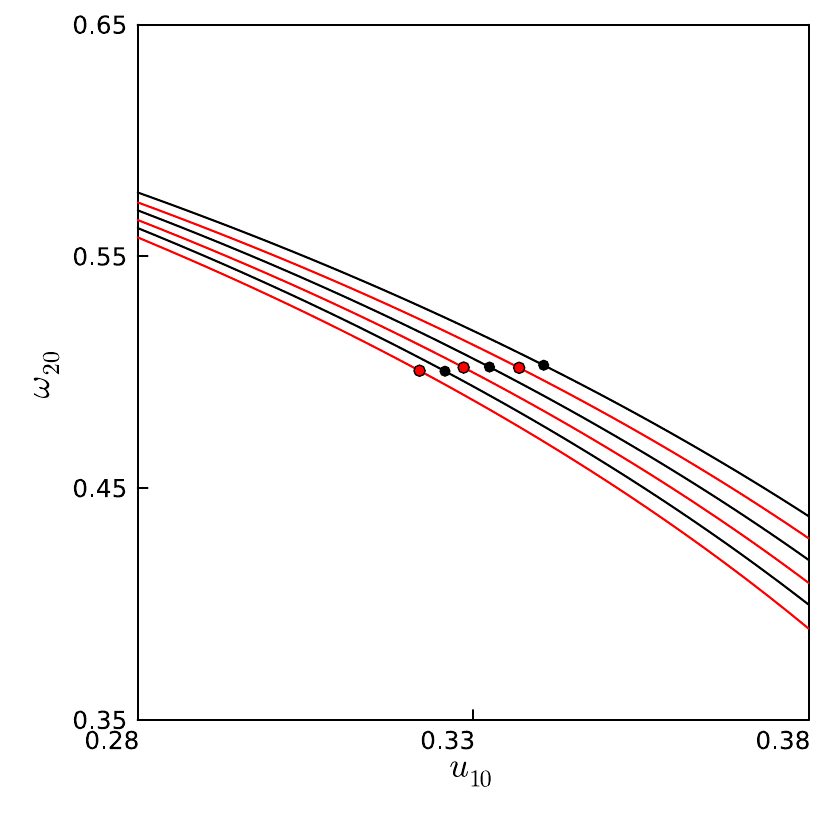}
\includegraphics[width=75mm]{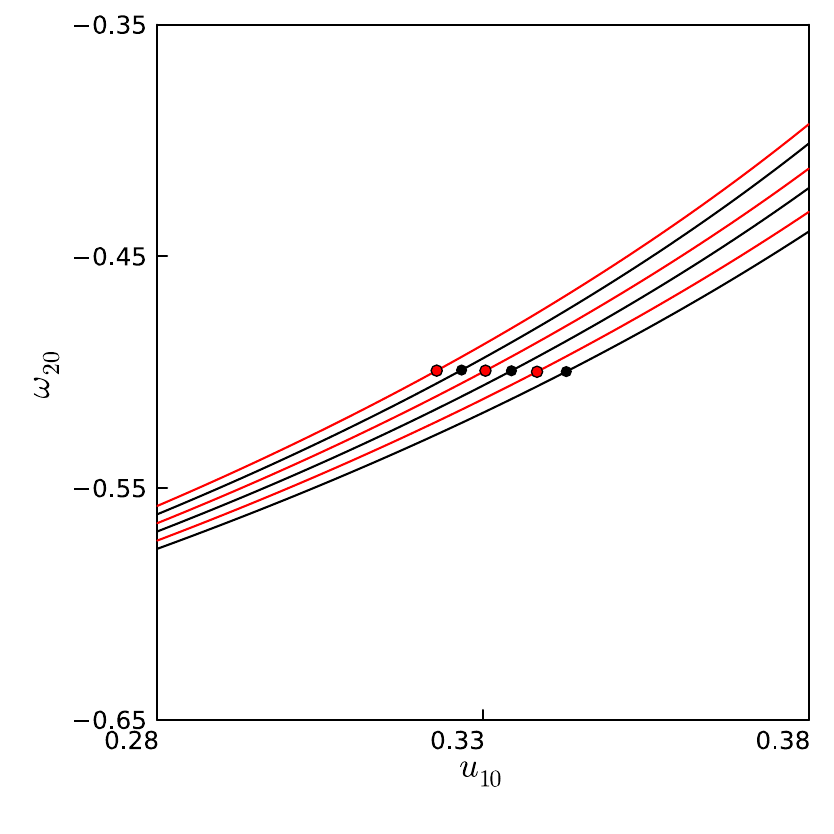}
\caption{Characteristic curves of six monoparametric families within $\widehat C_R$. In the second panel of the first appears from right to left families $1,\dots,6$ (black and red for odd and even trajectories, respectively). Filled dots indicate initial conditions of symmetric periodic orbits, corresponding initial conditions are given in Table 1.}
\label{figCR}
\end{figure}

We found for $\widehat C_R$, regular symmetric characteristic curves that resemble half ellipses that start and end at collision. These characteristic curves belong to six different monoparametric families which shrink as they approach to collision. We label the families according to the value that takes $u_{10}$ when $\omega_{20} = 0$; $i=1$ corresponds to the maximum $u_{10}$, $i=2$ for the second greatest value and so on until the six families are labelled (see Fig. \ref{figCR}).

Some properties of the different computed families can be explained thinking in bounded orbits near collision, we point out them in the following. Near collision the main force by the primaries over the test particle is due to the body $1$ since the small relative separation. In particular, if the path of the body $1$ does not change abruptly, for shorts times the bounded orbits of the test particle are approximately  elliptic, as happens in a two-body problem. Thus, at the initial time, the corresponding major and minor axes of such ellipses are parallel to $x$ and $y$ axes, and more than that, the test particle and the body $1$ appear at periapsis/apoapsis configuration; this property is expected for any initial condition ${\bf z}_0$ with characteristic time $T_0$ which gives rise to a motion near collision and satisfy BPV \eqref{bvp-3}. Something similar happens at $t=T_0$, the difference with respect to the configuration at $t=0$ is that in general, the major and minor axes of this new ellipse are not parallel to $x$ and $y$ axes. Reducing horizontal position and increasing vertical velocity of ${\bf z}_0$ appropriately we will be able to find an initial condition, let us call it ${\bf z}_1$, such that the BPV \eqref{bvp-3} is fulfilled and that the test particle upon reaching time $T_1$ near $T_0$ will make a half lap less than the orbit corresponding to ${\bf z}_0$. Thus, if the orbit defined by ${\bf z}_0$ has a periapsis/apoapsis configuration at $t = T_0$, then the orbit associated to ${\bf z}_1$ possess a apoapsis/periapsis configuration at $t = T_1$. This explains why, points belonging to neighbouring characteristic curves of $\widehat C_R$ define orbits qualitatively similar, which only differ by a half lap around the body $1$. Moreover, by means of considering the orbital period of the ellipse $T_{a} = 2 \pi a^{3/2}$, and supposing that for the initial condition ${\bf z}_0$, the required time to complete $n_l$ laps is $2 n_l \pi a_0^{3/2}$, then from ${\bf z}_0$ to ${\bf z}_1$ the semimajor axis $a_0$ is reduced to the one $a_1$ in such a way that the new orbit complete $n_l+1/2$ laps in the same time which is approximately $2 (n_l+1/2) \pi  a_1^{3/2}$. These expressions lead to the relation
$$
a_1 \approx a_0 \left(1 - \frac{1}{2(n_l+1/2)} \right)^{2/3},
$$
which serves as an estimate of how the space between the characteristic curves $\widehat C_R$ near collision are reducing.

The symmetry around $\omega_{20} = 0$ can be explained using similar arguments. Let us suppose that we have obtained initial conditions ${\bf z}_0$ with the properties mentioned in the previous paragraph. Notice that if we use \eqref{rv-rel} as guide to obtain numerically ${\bf z}_0$, then the corresponding orbit is prograde and changing only the sign of \eqref{rv-rel}, we could obtain numerically the ${\bf z}_2$ which leads to a retrograde orbit. The change of sign in the velocity in \eqref{rv-rel}, implies that the corresponding initial conditions to ${\bf z}_0$ and ${\bf z}_2$  in regularized coordinates only differ in a sign in the velocity, since the vectors \eqref{rv-rel} are relative to the body $1$. Thus, it is expected that the computed points within $\widehat C_R$ appear in quadruples, two in each representation, we mean, if $u_{10}$, $\omega_{20}$ ,$T_0$, $u_{20}=\omega_{10} = 0$ satisfies BPV \eqref{bvp-3}, also approximately points $u_{10},-\omega_{20},T_0$, $u_{20}=\omega_{10} = 0$ and $u_{20},\pm \omega_{10},T_0$, $u_{10}=\omega_{20} = 0$, satisfy the same BVP.

From the numerical results, we identify five regions that conform the computed characteristic curves within $\widehat C_R$. The regions $0<u_{10}<<1$ with $\omega_{20} >0$ and $\omega_{20} <0$ lead to prograde/retrograde orbits, respectively (see Fig. \ref{figCR}). On the other hand, the characteristic curves such that $\omega_{20} \approx 0$, describe the transition between prograde/retrograde orbits; in this case the corresponding initial conditions give rise to orbits with high eccentricity (see Fig. \ref{fig_high_e}). The other two regions are defined by very low values of $u_{10}$ such that the $|\omega_{20}|$ approaches to $1/\sqrt{2} \approx 0.707106781$ (with this, we have convergence for the binding energy \eqref{h-binding} as the curves that conform $\widehat C_R$ tend to collision). In the later case, the prograde/retrograde orbits computed with initial conditions within the characteristic curves that conform $\widehat C_R$ have a small eccentricity in the two-body problem approximation (see Fig. \ref{fig_low_e}).

\newpage

\subsubsection{Symmetric periodic orbits}

We have computed initial conditions of symmetric periodic orbits as those points $(u_{10},u_{20},\omega_{10},\omega_{20},\tau_0) \in \widehat C_{R}$ such that $\int_0^{\tau_0} {\bf u}^2 d\tau = 6 \overline{T}$ holds (see Section \ref{BVPregularized}).

\begin{figure}[h!]
\centering
\includegraphics[width=75mm]{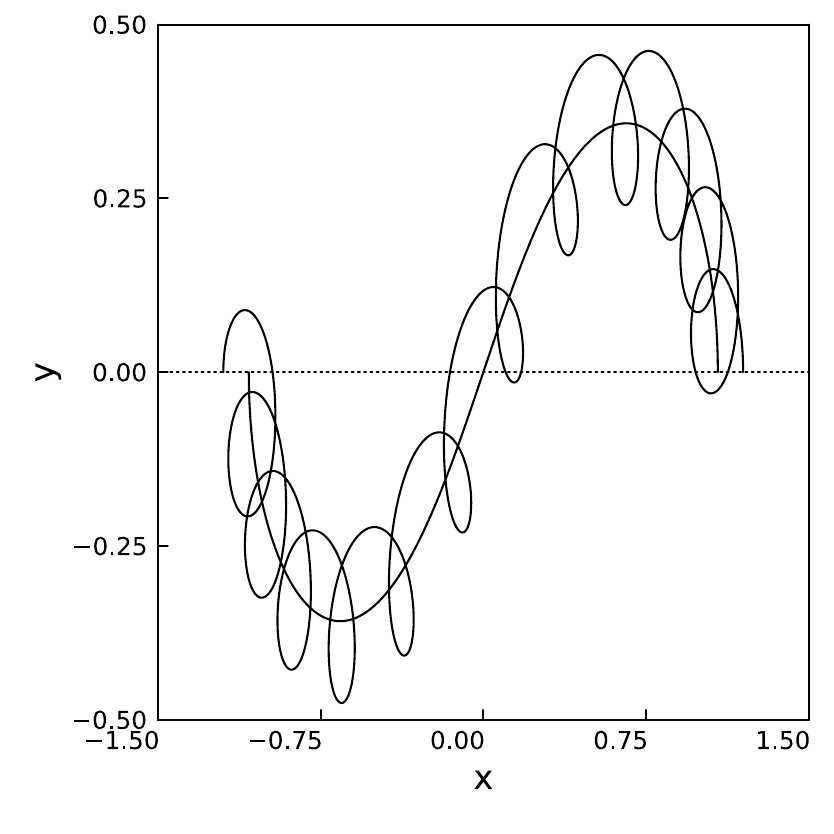}
\includegraphics[width=75mm]{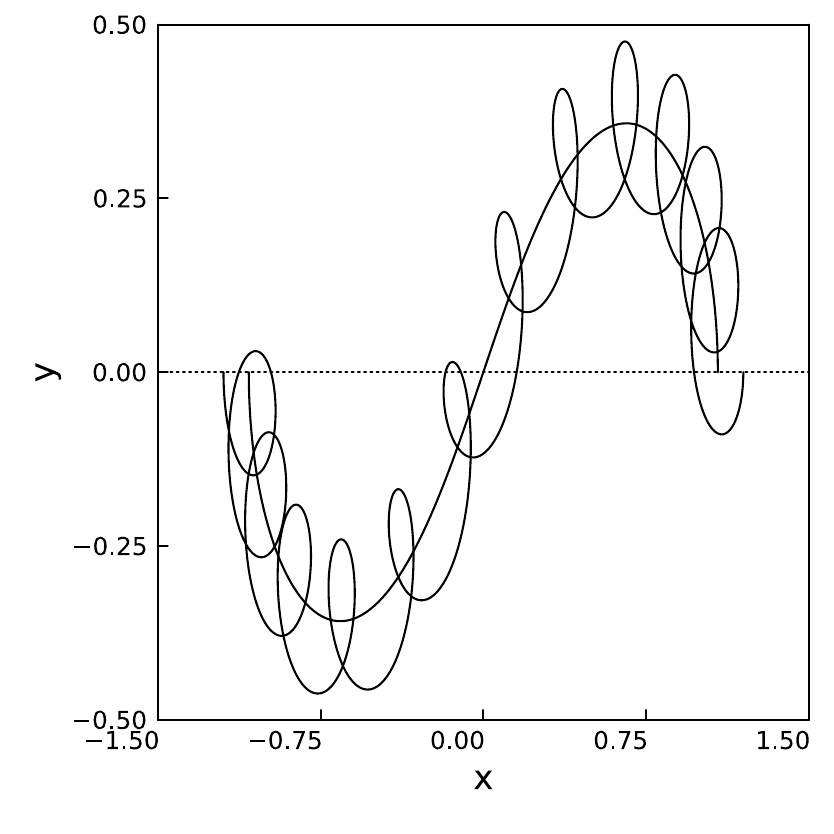}
\caption{Symmetric periodic orbits with low eccentricity - from left to right prograde and retrograde, respectively. The corresponding initial conditions are depicted in the bottom panels of Fig. \ref{figCR} for the largest values of $u_{10}$ (one value per panel).}
\label{fig_low_e}
\end{figure}

\begin{figure}[h!]
\centering
\includegraphics[width=75mm]{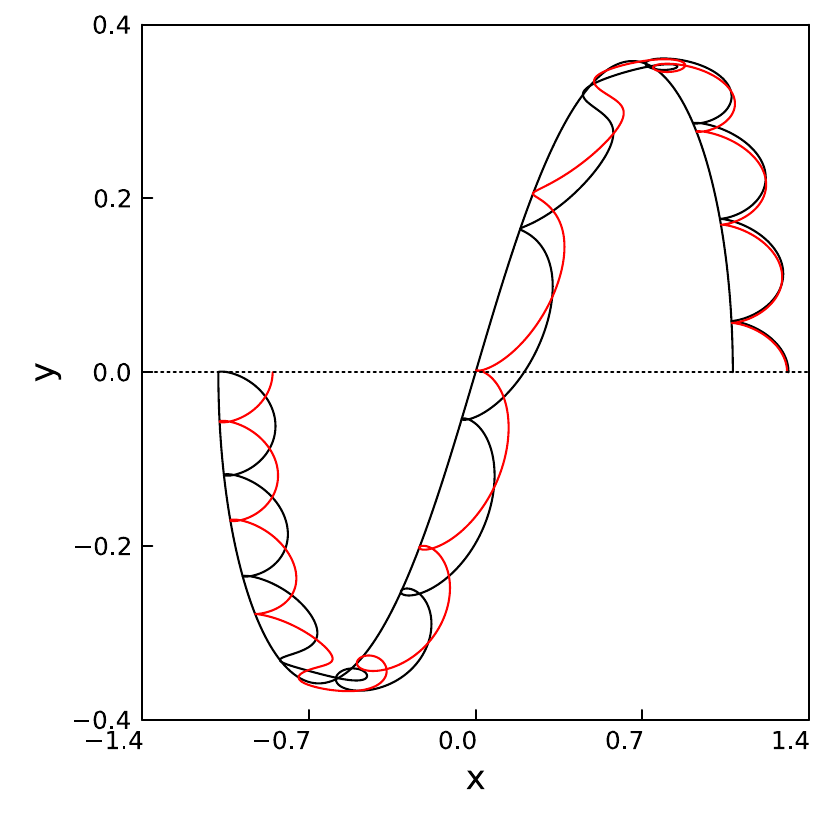}
\includegraphics[width=75mm]{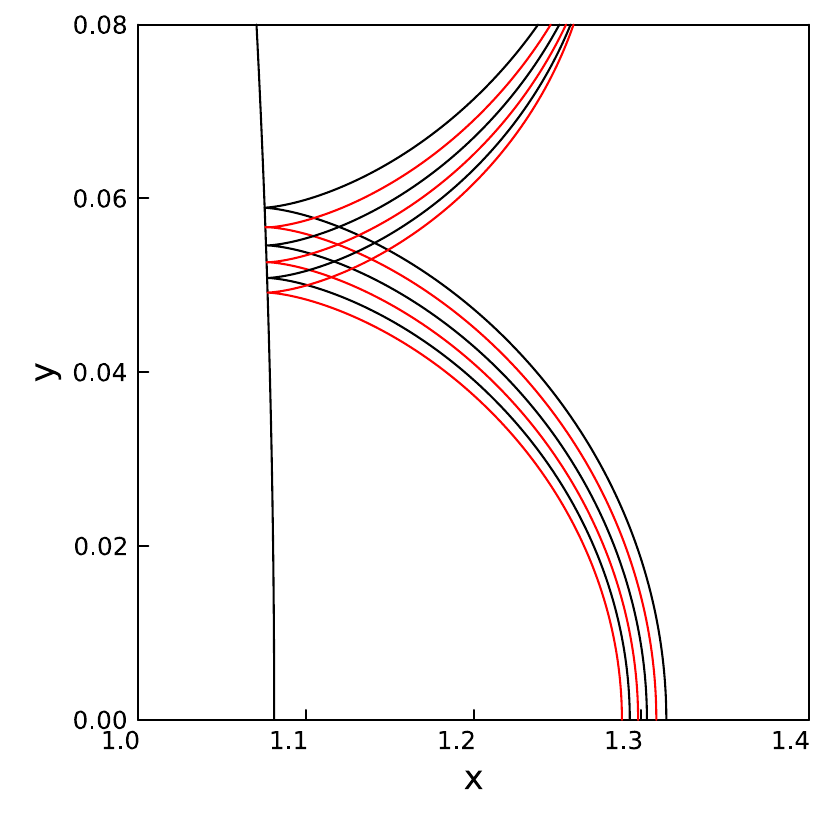}
\includegraphics[width=75mm]{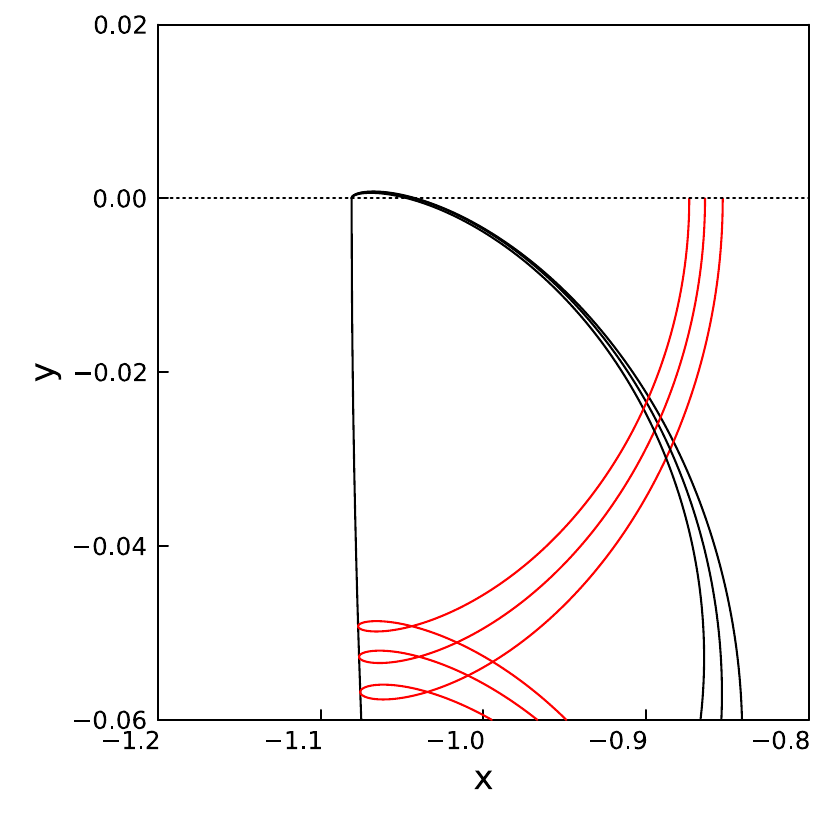}
\includegraphics[width=75mm]{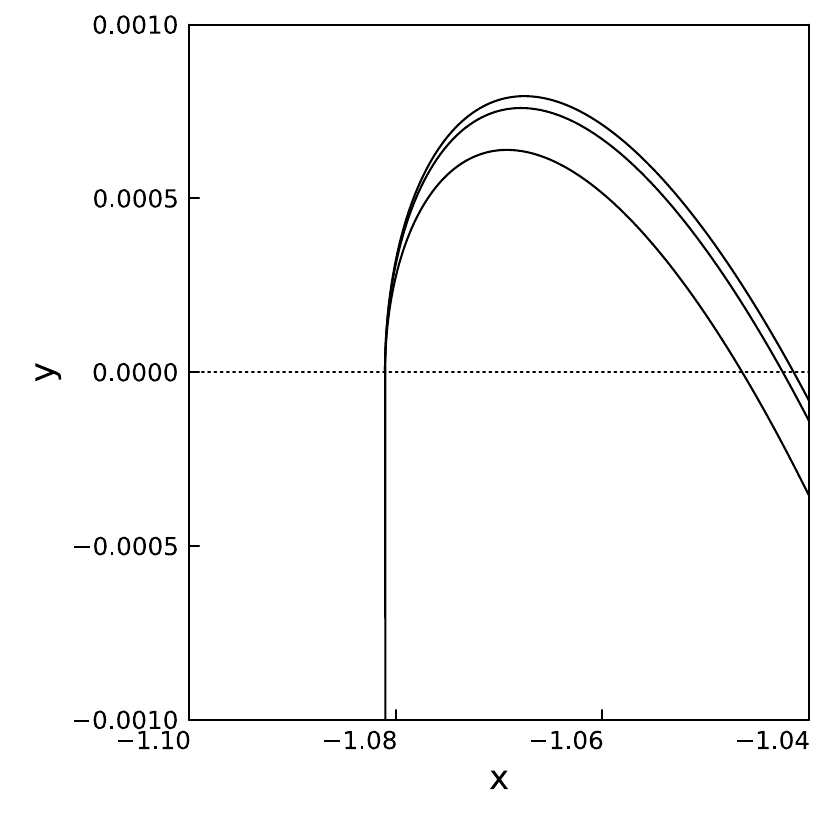}
\caption{Symmetric periodic orbits with high eccentricity. From top to bottom, left to right, in the first panel we shown two periodic orbits of the first row of Fig. \ref{figCR} - black trajectories for the biggest value of $u_{10}$ and red trajectories for the second largest value of $u_{10}$. The remaining three panels contain orbits whose initial condition appear in the upper second panel of Fig. \ref{figCR}. Numerical values of initial conditions are given in Table \ref{table1}.}
\label{fig_high_e}
\end{figure}

Each one of the six families present four points which define initial conditions of symmetric periodic orbits. These points appear in specific regions within $\widehat C_{R}$ as we explain in the following. In the region $0<u_{10}<<1$ with $\omega_{20} >0$, we have found prograde periodic orbits. In a similar way, for $0<u_{10}<<1$ with $\omega_{20} <0$, we have found retrograde periodic orbits. Example of these periodic orbits are shown in Fig. \ref{fig_low_e}.  On the other hand, for $\omega_{20} \approx 0$ only appears periodic prograde orbits, instead of both prograde and retrograde as happens in the case $0<u_{10}<<1$ (see Fig. \ref{fig_high_e}); this prevalence happens because for $\omega_{20} \approx 0$ prograde orbits are more stable than the retrogade ones (if the first body is fixed there is no prevalence of prograde or retrograde motion, but in the REFBP the first body is moving).

\begin{figure}[h!]
\centering
\includegraphics[width=75mm]{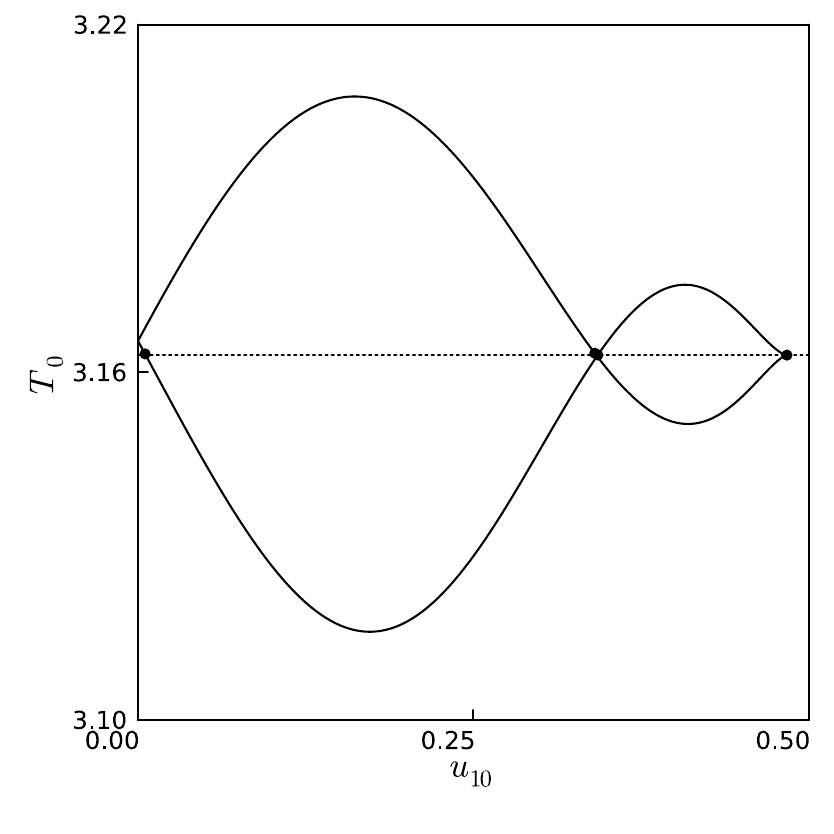}
\includegraphics[width=75mm]{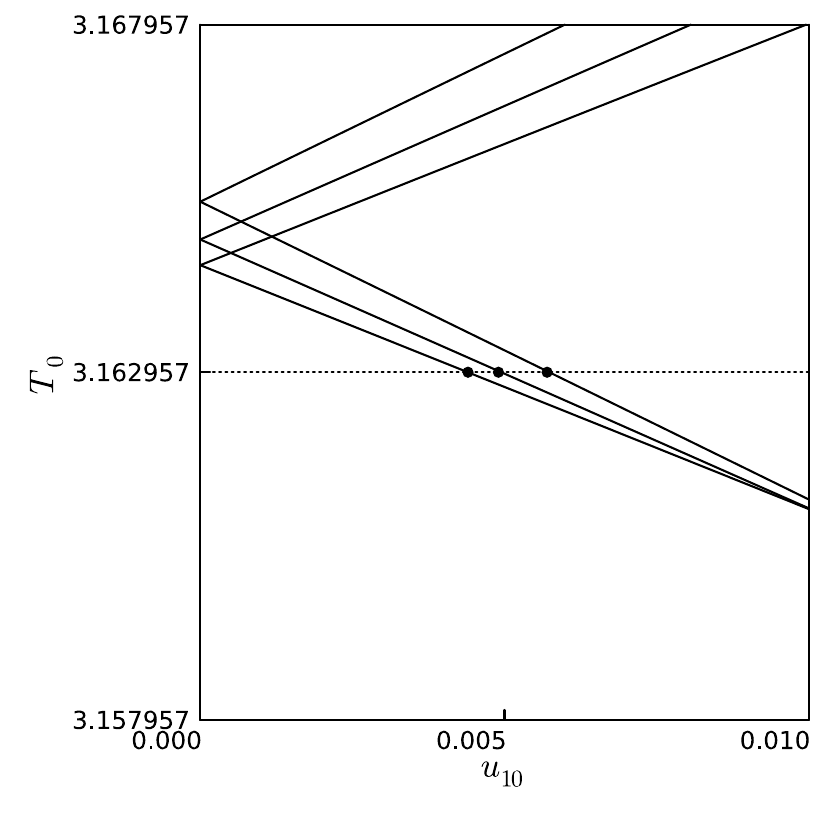}
\caption{Graphics of the characteristic time $T_0$ versus $u_{10}$. At left, first family within $\widehat C_R$, at right families 1, 3 and 5. Numerical values of initial conditions of periodic are given in Table 1.}
\label{figT1}
\end{figure}

\begin{figure}[h!]
\centering
\includegraphics[width=75mm]{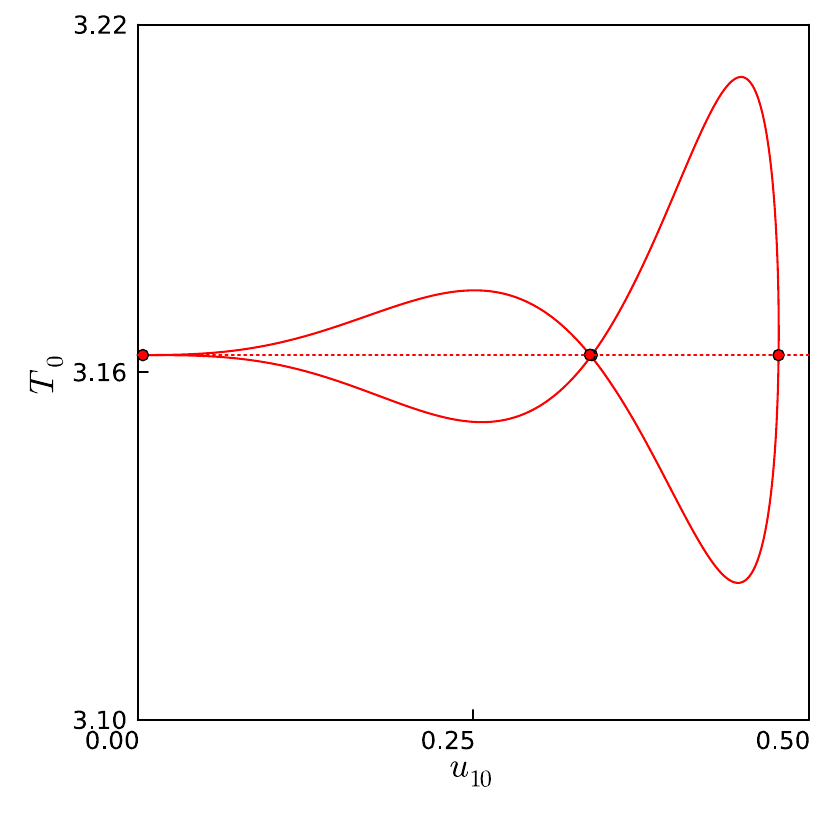}
\includegraphics[width=75mm]{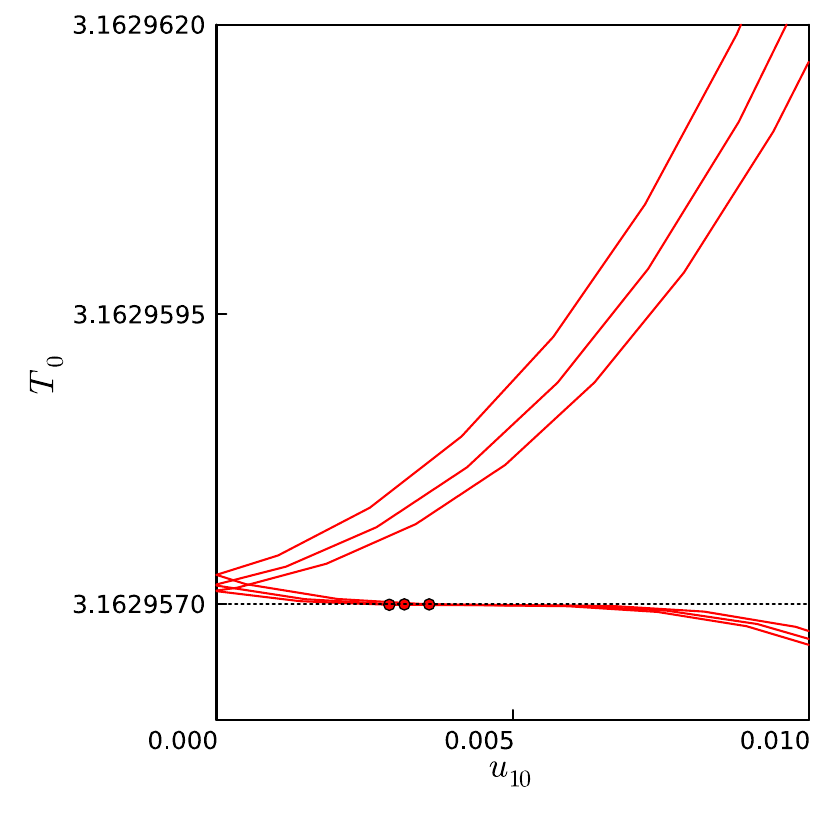}
\caption{Graphics of the characteristic time $T_0$ versus $u_{10}$. At left, second family within $\widehat C_R$, at right families 2, 4 and 6.  Numerical values of initial conditions of periodic are given in Table 1.}
\label{figT2}
\end{figure}

It only remains to describe the case when $u_{10} \to 0$. It is interesting to observe that in the retrograde case $u_{10} \to 0$,  $\omega_{20}>0$ does not appear initial conditions within $\widehat C_{R}$ that lead to periodic orbits. This fact can be appreciated in the graphics of the characteristic time $T_0$ versus initial position $u_{10}$ shown in Figs. \ref{figT1} and \ref{figT2} and Table \ref{table1}. The periodic orbits associated to initial points within $\widehat C_{R}$ with $u_{10} \to 0$, $\omega_{20} <0$, have the characteristic that the orbit of the test particle is close to collision for all time. Another interesting feature is that there are two shapes for $T_0$ as function of $u_{10}$, depending on the parity of the index of the family, as can be seen in Figs. \ref{figT1} and \ref{figT2}.

\begin{table}
\small
\begin{tabular}{ |p{0.25cm}|p{4.5cm}|p{4.5cm}|p{4.5cm}|  } \hline
\# &  $\tau_0$ &  $u_{10}$ & $\omega_{20}$ \\ \hline
1 & 2.695936701258381 $\times 10^{1}$ & 4.83740851834369 $\times 10^{-1}$ & 7.890047146086670  $\times 10^{-4}$ \\
2 & 2.767051193282881 $\times 10^{1}$ & 4.77531885645928 $\times 10^{-1}$ & 4.129816302813623  $\times 10^{-4}$ \\
3 & 2.837272976592164 $\times 10^{1}$ & 4.71631770925781 $\times 10^{-1}$ & 1.876247353294635  $\times 10^{-4}$ \\ 
4 & 2.906647017173428 $\times 10^{1}$ & 4.66013513823255 $\times 10^{-1}$ & 1.899368012266316  $\times 10^{-4}$ \\
5 & 2.975213137350447 $\times 10^{1}$ & 4.60650874021281 $\times 10^{-1}$ & 7.890042540965638  $\times 10^{-5}$ \\
6 & 3.043005316558202 $\times 10^{1}$ & 4.55526822164964 $\times 10^{-1}$ & -3.389574010239378 $\times 10^{-5}$ \\

7 & 2.693451258568199 $\times 10^{1}$ & 3.41577606593326 $\times 10^{-1}$ & 5.014586989085438 $\times 10^{-1}$ \\
8 & 2.764701267223720 $\times 10^{1}$ & 3.37068395755332 $\times 10^{-1}$ & 5.015850632582755 $\times 10^{-1}$ \\
9 & 2.835048393967071 $\times 10^{1}$ & 3.32904558862743 $\times 10^{-1}$ & 5.015290220699412 $\times 10^{-1}$ \\
10 & 2.904536554559997 $\times 10^{1}$ & 3.290176746761824 $\times 10^{-1}$ & 5.013584489391574 $\times 10^{-1}$ \\
11 & 2.973206096638524 $\times 10^{1}$ & 3.253033674527287 $\times 10^{-1}$ & 5.012030322717776 $\times 10^{-1}$ \\ 
12 & 3.041094179588408 $\times 10^{1}$ & 3.217154560535076 $\times 10^{-1}$ & 5.011123753102673 $\times 10^{-1}$ \\

13 & 2.693451258568199 $\times 10^{1}$ & 3.427121370863340 $\times 10^{-1}$ &-4.998145423262709 $\times 10^{-1}$ \\
14 & 2.764701322142713 $\times 10^{1}$ & 3.382737664718438 $\times 10^{-1}$ &-4.998149628532847 $\times 10^{-1}$ \\
15 & 2.835048393967071 $\times 10^{1}$ & 3.342059257610913 $\times 10^{-1}$ &-4.995921597981327 $\times 10^{-1}$ \\
16 & 2.904536529524823 $\times 10^{1}$ & 3.300847221928169 $\times 10^{-1}$ &-4.997507685924308 $\times 10^{-1}$ \\
17 & 2.973206096638524 $\times 10^{1}$ & 3.262535137756929 $\times 10^{-1}$ &-4.997539033212443 $\times 10^{-1}$ \\ 
18 & 3.041094181647387 $\times 10^{1}$ & 3.226323426319703 $\times 10^{-1}$ &-4.996974241102303 $\times 10^{-1}$ \\

19 & 2.695936701257934 $\times 10^{1}$ & 5.735296308245588 $\times 10^{-3}$ & -7.070568089673658 $\times 10^{-1}$ \\
20 & 2.767051919036603 $\times 10^{1}$ & 5.138142721366328 $\times 10^{-3}$ & -7.070656442349244 $\times 10^{-1}$ \\
21 & 2.837272976592157 $\times 10^{1}$ & 4.945918740892248 $\times 10^{-3}$ & -7.070677206777947 $\times 10^{-1}$ \\
22 & 2.906646918450844 $\times 10^{1}$ & 4.550398110681072 $\times 10^{-3}$ & -7.070729270229507 $\times 10^{-1}$ \\
23 & 2.975213137350443 $\times 10^{1}$ & 4.396529780308043 $\times 10^{-3}$ & -7.070744473086668 $\times 10^{-1}$ \\
24 & 3.043005621836102 $\times 10^{1}$ & 4.423458624321410 $\times 10^{-3}$ & -7.070733177487463 $\times 10^{-1}$ \\
\hline
\end{tabular}
\caption{Initial conditions of symmetric $12\overline{T}$ periodic orbits in the REFBP.}
\label{table1}
\end{table}

\newpage

\section{Conclusions}
We have studied periodic motions near collision in a restricted problem which consists of a massless particle moving under the gravitational influence due to three bodies following the figure-eight choreography.

It is well known that the periodic orbits of non-autonomous systems are isolated since, in a generic way, non-autonomous systems do not possess a first integral integral, unlike autonomous cases, for example the circular restricted three-body problem. Our approach is to use boundary problems whose solutions are conformed by monoparametric families in which appear initial conditions of symmetric periodic orbits.

We found two types of orbits near collision. In the moon case the numerical results suggest the existence of a succession of orbits which increase their number of laps around the primary 
as the distance between the test particle and the primary decreases. This result is consistent with \cite{Ortega} where the authors found regularized periodic solutions of the test particle in a neighbourhood of one primary.

In our study, the regularization of the equations of motion was important for the numerical aspect. Although using regularized coordinates gives rise to a regular system, there are significant numerical changes in orbits near collision. Moreover, the system is non-autonomous and we do not have an explicit parametrization for  state vectors of the primaries, therefore the computation of such state vectors involves both natural and fictitious time variables which leads to another numerical complication.




\section*{Acknowledgements}
First and third authors are pleased to acknowledge Asociación Mexicana de Cultura A. C. and National System of Researchers (SNI).


\begin{thebibliography}{9}

\bibitem{BarreraBengocheaGargia} Barrera, C., Bengochea, A., Garc\'ia-Azpeitia, Comet and Moon Solutions in the Time--Dependent Restricted $(n+1)$--Body Problem. Journal of Dynamics and Differential Equations. \textbf{34}, 1187-1207 (2022)

\bibitem{Broucke} Broucke R. On Changes of Independent Variable in Dynamical 		Systems and Applications to Regularization.  Icarus. \textbf{7}, 221–231 (1967)

\bibitem{BurgosVI} Burgos-Garc\'ia J., Bengochea A., Horseshoe orbits in the restricted four body problem. Astrophys Space Sci \textbf{362}, 212-226 (2017) 


\bibitem{BurgosBengocheaPerez} Burgos-Garc\'ia, J., Bengochea, A.,  Franco-P\'erez, L. The Spatial Hill Four--Body Problem I - An 
Exploration of Basic Invariant Sets.  Communications in Nonlinear Science and Numerical Simulation. \textbf{108}. (2022)





\bibitem{BurgosDelgado} Burgos-Garc\'ia, J., Delgado, J.; Periodic orbits in the restricted four-body problem with two equal masses. Astrophys Space Sci \textbf{345}, 247–263 (2013). https://doi.org/10.1007/s10509-012-1118-2

\bibitem{BurgosDelgadoII} Burgos-Garc\'ia, J., Delgado, J.; On the ``Blue sky
catastrophe'' termination in the restricted four body problem. Celest Mech Dyn Astr \textbf{117}, 113-136 (2013). https://doi.org/10.1007/s10569-013-9498-3.

\bibitem{BurgosLessardJames} Burgos-Garc\'ia, J., Lessard, J.,  James, J.D.M. Spatial periodic orbits in the equilateral circular restricted four-body problem: computer-assisted proofs of existence. Celest Mech Dyn Astr \textbf{131}, 2 (2019). https://doi.org/10.1007/s10569-018-9879-8




\bibitem{Celletti} Celletti, A. Stability and Chaos in Celestial Mechanics.  Springer. (2010)

\bibitem{Chenciner-Montgomery} Chenciner, A., Montgomery , R. A remarkable periodic solution of the three-body problem in the case of equal masses. Annals of Mathematics. \textbf{152}, 881-901 (2000)

\bibitem{Fujiwara} Fujiwara T., Fukuda H. and Ozaki H. Choreographic three bodies on the lemniscate.  J. Phys. A: Math. Gen.n. \textbf{36}, 2791–2800 (2003)

\bibitem{Lamb} Lamb, S. W., Roberts, J. A. G., Time-reversal symmetry in dynamical systems: A survey. Physica D. \textbf{112}, Issues 1-2, 1-39 (1998)

\bibitem{LaraBengochea} Lara, R., Bengochea, A., A Restricted Four-body Problem for the Figure-eight Choreography. Regular and Chaotic Dynamics. \textbf{26}, 222-235 (2022)

\bibitem{Leandro} Leandro, E.S.G.: On the central configurations of the planar restricted four-body problem. J. Differ. Equ. \textbf{226(1)}, 323–351 (2006)

\bibitem{KepleyII} Kepley S, Mireles James, J.D.: Homoclinic dynamics in a restricted four-body problem: transverse connections for the saddle-focus equilibrium solution set. Celest Mech Dyn Astr \textbf{131}, 13 (2019).

\bibitem{MunozFreireGalanVader} Mu\~noz-Almaraz, F.J., Freire, E., Galán, J., Vanderbauwhede, A., Continuation of normal doubly symmetric orbits in conservative reversible systems. Celestial Mechanics and Dynamical Astronomy. \textbf{97}, 17-47 (2007)


\bibitem{Ortega} Ortega, R., Zhao, L., Generalized periodic orbits in some restricted three-body problems. Z. Angew. Math. Phys. \textbf{72}, 40 (2021)


\bibitem{Szebehely} Bettis, D.G., Szebehely, V. Treatment of close approaches in the numerical integration of the gravitational problem of N bodies.  Astrophysics and Space Science. \textbf{14}. (1971)

\bibitem{Moore} Moore, C. Braids in Classical Gravity. Phys. Rev. Lett. \textbf{70}, 3675–3679 (1993)

\bibitem{PDS1} Balsalobre-Ruza O, et al. Tentative co-orbital submillimeter emission within the Lagrangian region L5 of the protoplanet PDS70b.  Astronomy \& Astrophysics \textbf{675}, 1–8 (2023)

\bibitem{alma} Does this exoplanet have a sibling sharing the same orbit? Atacama Large Milimiter Array (ALMA). Press Releases. (2023) \url{https://www.almaobservatory.org/en/press-releases/does-this-exoplanet-have-a-sibling-sharing-the-same-orbit/?fbclid=IwAR26t87nIvnCzKpDzF6RsbvKt4H6XBVnr22hSGdXOuZIGblBeAsDLtfphAs}




\end{thebibliography}
\end{document}